\def\g{{\gamma}}

\def\xi{{\xi}}
\def\p{{\pi}}
\def\r{{\rho}}
\def\s{{\sigma}}
\def\t{{\tau}}

\def\De{{\Delta}}
\def\Th{{\Theta}}
\def\Xi{{\Xi}}
\def\Pi{{\Pi}}

\def\F{{\cal F}}

\def\ra{{\rightarrow}}

\documentclass[10pt]{article}
\usepackage{psfig,latexsym,amsbsy}
\newtheorem{thm}{Theorem}
\newtheorem{cor}[thm]{Corollary}

\newtheorem{quest}[thm]{Question}

\begin{document}

\baselineskip=24pt
%(makes the paper double-spaced)

% #######################################################################
% #######################################################################
%	
%	TITLE PAGE: 
%	
\title{Cover Pebbling Hypercubes}
\author{
 \\
 \\
Glenn H. Hurlbert\thanks{Partially supported by National Security
        grant \#MDA9040210095.}\\
Department of Mathematics and Statistics\\
Arizona State University\\
Tempe, AZ 85287-1804\\
email: hurlbert@asu.edu\\
 \\
and \\
 \\
Benjamin Munyan\\
Department of Mathematics and Statistics\\
Arizona State University\\
Tempe, AZ 85287-1804\\
email: munyan@asu.edu\\
}
\maketitle
\newpage

%\small\normalsize
%\renewcommand{\baselinestretch}{2.5}
%(makes the paper double-spaced)

% #######################################################################
% #######################################################################
%	
%	ABSTRACT:
%	
\begin{abstract}
Given a graph $G$ and a configuration $C$ of pebbles on the vertices
of $G$, a pebbling step removes two pebbles from one vertex and places
one pebble on an adjacent vertex.
The cover pebbling number $\g=\g(G)$ is the minimum number so that every 
configuration of $\g$ pebbles has the property that, after some sequence 
of pebbling steps, every vertex has a pebble on it.
We prove that the cover pebbling number of the $d$-dimensional
hypercube $Q^d$ equals $3^d$.
\vspace{0.2 in}

\vspace{0.2 in}

\end{abstract}
\newpage

% #######################################################################
% #######################################################################
%	
%	BEGINNING OF PAPER:
%	
\section{Introduction}\label{intro}

Given a graph $G$ and a configuration $C$ of pebbles on the vertices
of $G$, a pebbling step removes two pebbles from one vertex and places
one pebble on an adjacent vertex.
The {\it pebbling number} $\p=\p(G)$ is the minimum number so that every
configuration of $\p$ pebbles has the property that, for any given
{\it target} vertex, after some sequence of pebbling steps, the target
has a pebble on it.
Two basic results are that the complete graph satisfies $\p(K_n)=n$, 
and that the path satisfies $\p(P_n)=2^{n-1}$, where $n=n(G)$ is the number
of vertices of the graph $G$.
The pebbling numbers of trees and cycles have also been computed in
\cite{M} and \cite{PSV}, respectively.
It is easy to see that $n(G)$ and $2^{{\rm diam}(G)}$ are each lower 
bounds on $\p(G)$, and $K_n$ and $P_n$ show that these bounds can be tight
(here ${\rm diam}(G)$ is the diameter of $G$).
Chung \cite{C} proved that the $d$-dimensional hypercube, or $d$-cube,
satisfies $\p(Q^d)$, which interestingly is tight on both accounts.
There is a radiply growing literature on the subject (see \cite{H}),
including a handful of variations on the theme such as optimal pebbling and
pebbling thresholds.

In this paper we consider the {\it cover pebbling number}, first introduced
in \cite{CCF}.
The cover pebbling number $\g=\g(G)$ is the minimum number so that every 
configuration of $\g$ pebbles has the property that, after some sequence 
of pebbling steps, every vertex has a pebble on it.
Crull, et al. \cite{CCF}, find the cover pebbling number of trees and
complete graphs.
Because $n$ targets must be reached instead of just one, the bound
$\g(G)\le n\p(G)$ holds in general.
In light of Chung's result this yields $\g(Q^d)\le 4^d$.
Here we prove the following.

\begin{thm}\label{cube}
The cover pebbling number of the $d$-cube is $\g(Q^d)=3^d$.
\end{thm}

In \cite{CCF} is also defined the {\it cover pebbling ratio} $\r=\g/\p$.
Using Moews's result on trees, the authors show that $\r$ can be as small 
as 2 (cliques, paths) and as large as $n/\lg n$ (brooms or fuses).
In the case of cubes we have the following.

\begin{cor}\label{ratio}
The cover pebbling ratio of the $d$-cube is 
$\r(Q^d)=n^{\lg 3-1}=n^{.5849625\ldots}$.
\end{cor}

% #######################################################################
% #######################################################################
%	
%	PRELIMINARIES:
%
\section{Preliminaries}\label{prelim}

We begin by developing the terminology we will employ.
A {\it configuration} $C$ of pebbles on the vertices $V(G)$ of a graph 
$G$ is a function $C:V(G)\ra \{0,1,2,\ldots\}$, where $C(x)$ is the 
number of pebbles on vertex $x$.
The {\it size} of $C$ equals $|C|=\sum_{x\in V}C(x)$, the total number 
of pebbles on $G$.
The {\it support} of $C$ is the set $\s=\s(C)=\{x\in V\ |\ C(x)>0\}$
of vertices having at least one pebble.
We say that $C$ is {\it simple} if $|\s(C)|=1$, is a {\it cover} if
$\s(C)=V(G)$, and is {\it even} if $C(x)$
is even for every $x$.
A vertex $x$ is {\it empty} if $C(x)=0$, a {\it one} if $C(x)=1$, a 
{\it two} if $C(x)=2$, and is {\it large} if $C(x)\ge 2$.
The configuration $C$ is called {\it coverable} if after some sequence of 
pebbling steps no vertex is empty.
If $C$ has a large vertex having an empty neighbor, we say that $C$ is
{\it open}; otherwise it is {\it closed}.
Finally, for the purposes of this article, $C$ is {\it good} if 
$|C|\ge 3^n-|\s(C)|+1$, and is {\it sharp} if equality holds.

We note that a simple configuration of size $3^n-1$ is not coverable.
Indeed, to reach a vertex at distance $i$ from a simple support requires 
$2^i$ pebbles, as shown by $\p(P_{i+1})$.
Thus $\g(Q^d)\ge\sum_{i=0}^d {d\choose i}2^i=3^d$.

Our proof of Theorem \ref{cube} borrows an idea from Chung's proof that
$\p(Q^d)=2^d$ (see \cite{C}).
Using induction, she proved the extra statement that $Q^d$ had the
{\it 2-pebbling property}.
A graph $G$ has this property if, from every configuration $C$ that
satisfies $|C|\ge 2\p(G)-|\s(C)|+1$, one can place 2 pebbles on any
specified target.
It is a curious property that suggests that more concentrated
configurations require greater size to maintain power.
For example, if $C$ is a cover on $Q^d$ then it needs only to be of size
$n+1$ to 2-pebble an arbitrarily chosen target, while if it is simple 
then it needs to be of size $2n$ instead.
Somewhat analogously we prove the following.

\begin{thm}\label{good}
Every good configuration on $Q^d$ is coverable.
\end{thm}

It is clear that Theorem \ref{good} implies Theorem \ref{cube}.
It is also true that we need only prove Theorem \ref{good} in the case
that $C$ is sharp.

% #######################################################################
% #######################################################################
%	
%	PROOF OF THEOREM:
%
\section{Proof of Theorem \ref{good}}\label{proof}

As noted, we may assume that the configuration $C$ is sharp.
If $C$ is open then the appropriate pebbling step creates a
sharp configuration of smaller size.
In this case we use induction on sharp configuration size.
The base case has support equal to the cube: it is a cover.
Thus we may assume that the configuration is closed.

If there is only one large vertex then we use induction on support size.
The base case has support size $d+1$ and is easily coverable: only 
$1+\sum_{i=2}^d {d\choose i}2^i=3^d-2d$ pebbles are needed on the large
vertex and there are exactly $(3^d-(d+1)+1)-d=3^d-2d$ pebbles on it.
In general, we compare the configuration to the one obtained by removing 
a one and placing two pebbles on the large vertex.
Because the new configuration is also sharp and has smaller support size 
it is coverable.
Since it takes at least two pebbles to cover the newly emptied vertex,
the remaining pebbles can cover all other empties.
Thus those pebbles can still cover all other empties in the original
configuration --- that is, $C$ is coverable.
Thus we may assume there are at least two large vertices.

Because the configuration is closed this means that the support size is
at least $|N(u)|+|N(v)|-|N(u)\cap N(v)|\ge 2(n+1)-2=2n$,
where $N(x)$ is the closed neighborhood of $x$.
Hence, for $d\le 2$ the configuration is already a cover.
Now assume $d=3$.

If there are two antipodal large vertices then $C$ is a cover.
Otherwise if there are two adjacent large vertices then $|\s(C)|\ge 6$.
We will assume that there are just two large vertices in this case;
it is simpler to make the arguments with more large vertices.
Suppose $|\s(C)|=6$, then $|C|=3^3-6+1=22$.
If the largest vertex has at least 13 pebbles then it alone can cover 
the two empties, so let's say it has at most 12 pebbles.
This means that the other large has at least 6 pebbles, meaning it can
cover its closest empty vertex.
Since the largest vertex must also have at least 6 pebbles, it can cover
the remaining empty.
It is even simpler to argue that such a configuration of support size 7
is coverable.
Now we may assume there are no adjacent large vertices.

For the case of just two large vertices, one can argue along similar lines 
as above that configurations of support size 6 or 7 are coverable.
For the case of three large vertices the arguments are easier.
A configuration having four large vertices, none of which are adjacent 
or antipodal, is a cover.
Now we may assume that $d\ge 4$, and we will argue by induction on $d$.

We use the natural labelling of $V(Q^d)$ by binary $d$-tuples, with
adjacent vertices determined by Hamming distance 1.
For any coordinate $j$ we can {\it cut} $Q^d$ into two copies of $Q^{d-1}$,
the {\it top} copy $T=T_j$ having $j^{\rm th}$ coordinate equal to 1 and
the {\it bottom} copy $B=B_j$ having $j^{\rm th}$ coordinate equal to 0.

For a given top/bottom cut of the cube, if both corresponding configurations
$C_T$ and $C_B$ are good then we are done by induction.
Thus we may assume that every such cut has $|C_B|=3^{n-1}-|\s_B|+1-\De$
for some $\De>0$
(by swapping 0 and 1 on a cut (coordinate) if necessary, we can make sure
that it is the bottom rather than the top that fails the condition).
Hence 
\begin{eqnarray*}
|C_T|&=&(3^n-|\s|+1)-|C_B|\\
&=&(3^n-|\s|+1)-(3^{n-1}-|\s_B|+1-\De)\\
&=&(3^{n-1}-|\s_T|+1)+(3^{n-1}+\De-1)\ .\\
\end{eqnarray*}

Since $n\ge 4$ we have $3^{n-2}\ge 2^{n-1}\ge |\s_T|$.
Thus $3^{n-1}\ge 3|\s_T|$ and so $3^{n-1}-|\s_T|+1\ge 2|\s_T|$.
Consider a the configuration $R_T$ of empties, ones and twos on $T$ 
that is congruent mod 2 to $C_T$ and has the same support.
It has size at most $2|\s_T|\le 3^{n-1}-|\s_T|+1$ and is domininated by $C_T$
(that is, $R_T(x)\le C_T(x)$ for all $x\in T$).
Thus the configuration $S_T=C_T-R_T$ is even and has size at least 
$3^{n-1}+\De-1$.
Because $C$ has at least two large vertices we have $|\s_B|\ge 2$.
Therefore $0\le |C_B|\le 3^{n-1}-1-\De$, and so $\De\le 3^{n-1}-1$.
This means that $3^{n-1}+\De-1\ge 2\De$ and so $\De$ pebbles from $S_T$
can be moved from $T$ to $B$.
This results in a good configuration on $T$ and a sharp configuration on
$B$, and the proof is finished by induction.

% #######################################################################
% #######################################################################
%	
%	OPEN QUESTIONS:
%	
\section{Open Questions}\label{open}

We have seen in the case of cubes that there is a simple non-coverable 
configuration of size $\g(Q^d)-1$.
This follows the behavior of complete graphs and cubes, and so we
reiterate a question first raised in \cite{CCF}.

\begin{quest}
Is it true that every graph $G$ has a noncoverable configuration of size
$\g(G)-1$ that is simple?
\end{quest}
If this is true then every graph $G$ would have cover pebbling number
$\g(G)=\max_{v\in V(G)}\sum_{x\in V(G)}2^{{\rm dist}(v,x)}$, where
${\rm dist}(v,x)$ is the distance between $v$ and $x$.

Define the cover pebbling ratio of a class $\F$ of graphs as 
$\r(\F)=\sup_{G\in\F}\r(G)$.
We noted that the class of complete graphs and paths has cover pebbling
ratio 2, and that the class of trees has cover pebbling ratio $n/\lg n$.
Here we discovered that the class of cubes has cover pebbling ratio
$n^{.58\ldots}$.
\begin{quest}
Is there an infinite class of graphs whose cover pebbling ratio is either
smaller than 2 or larger than $n/\lg n$?
\end{quest}

Another interesting pursuit is the following.
\begin{quest}
Is there an appropriate graph invariant that identifies (rather than
characterizes) either large or small cover pebbling ratio?
\end{quest}
As evidenced by paths, diameter is not such an invariant.

Graham's nototious pebbling conjecture states that 
$\p(G\Box H)\le\p(G)\p(H)$ for every pair of graphs $G$ and $H$, 
where $\Box$ is the cartesian product.
\begin{quest}
Is it true that every pair of graphs $G$ and $H$ satisfy
$\g(G\Box H)\le \g(G)\g(H)$?
\end{quest}
This is true for cubes, and because of the lack of need in the cubes for a 
special property like 2-Pebbling, this may be a simpler question to resolve.

For those who like probabilistic questions see, for example, \cite{BBCH}
for the definition of the pebbling threshold for a sequence of graphs.
The cover pebbling threshold is defined analogously.
\begin{quest}
Is the cover pebbling threshold for a graph sequence equal to $\Th(\r\t)$?
\end{quest}

% #######################################################################
% #######################################################################
%	
%	BIBLIOGRAPHY:
%	
\bibliographystyle{plain}

\begin{thebibliography}{99}
\bibitem{BBCH} 
A. Bekmetjev, G. Brightwell, A. Czygrinow and G. Hurlbert,
{\it Thresholds for families of multisets, with an application to 
graph pebbling},
Discrete Math. {\bf 269} (2003), no.1-3, 21--34.
\bibitem{C} 
F.R.K. Chung
{\it Pebbling in hypercubes},
SIAM J. Disc. Math. {\bf 2} (1989), 467--472.
\bibitem{H} 
G. Hurlbert,
{\it A survey of graph pebbling},
Congr. Numer. {\bf 139} (1999), 41--64.
(Proceedings of the Thirtieth Southeastern International Conference on 
Combinatorics, Graph Theory, and Computing, Boca Raton, FL, 1999.)
\bibitem{CCF} 
B. Crull, T. Cundiff, P. Feldman, G. Hurlbert, L. Pudwell, S. Szaniszlo 
and Z. Tuza,
{\it The cover pebbling number of graphs},
submitted.
\bibitem{M} 
D. Moews,
{\it Pebbling graphs},
J. Combin. Theory (Ser. B) {\bf 55} (1992), 244--252.
\bibitem{PSV} 
L. Pachter, H. Snevily and B. Voxman,
{\it On pebbling graphs},
Congr. Numer. {\bf 107} (1995), 65--80.
\end{thebibliography}
%   

% #######################################################################
% #######################################################################
%	
%	END OF PAPER
%	
\end{document}